\documentclass{amsart}
\usepackage{color}
\usepackage[all]{xy}
\usepackage{amssymb}
\let\blb\mathbb

\def \ZZ{{\blb Z}}

\def \HH{{\blb H}}

\def\Hom{\operatorname {Hom}}

\def\r{\rightarrow}

\DeclareMathOperator{\Aut}{Aut}

\newtheorem*{question1}{Question}
\newtheorem*{example1}{Example}
\newtheorem*{lemma1}{Lemma 1}
\newtheorem*{lemma2}{Lemma 2}

\theoremstyle{definition}

{

}

\theoremstyle{remark}

\newtheorem*{remark1}{Remark}

\newdimen\uboxsep \uboxsep=1ex
\def\uboxn#1{\vtop to 0pt{\hrule height 0pt depth 0pt\vskip\uboxsep
\hbox to 0pt{\hss #1\hss}\vss}}

\def\uboxs#1{\vbox to 0pt{\vss\hbox to 0pt{\hss #1\hss}
\vskip\uboxsep\hrule height 0pt depth 0pt}}

\def\HH{\operatorname{HH}}
\def\Tw{\operatorname{Tw}}

\def\Perf{\operatorname{Perf}}
\let\oldmarginpar\marginpar
\def\marginpar#1{\oldmarginpar{\tiny #1}}
\def\Tw{\operatorname{tw}}
\title{A note on non-unique enhancements}
\author{Alice Rizzardo}
\email[Alice Rizzardo]{alice.rizzardo@ed.ac.uk}
\address{School of Mathematics\\The University of Edinburgh\\James Clerk Maxwell Building\\The King's Buildings\\Peter Guthrie Tait Road\\Edinburgh, EH9 3FD\\Scotland, UK}
\author{Michel Van den Bergh}
\email[Michel Van den Bergh]{michel.vandenbergh@uhasselt.be}
\address{Universiteit Hasselt\\ Universitaire Campus\\ 3590 Diepenbeek\\Belgium}
\thanks{The first author is a Postdoctoral Research Fellow at the University of Edinburgh. She is supported by EPSRC grant EP/N021649/1. The second author is a senior researcher at the Research Foundation - Flanders (FWO). He is supported by the FWO-grant G0D8616N ``Hochschild cohomology and deformation theory of triangulated categories.''}
\begin{document}
\begin{abstract}
We give an easy example of a triangulated category, linear over a field~$k$, with two different
enhancements, linear over $k$, answering a question of Canonaco and Stellari.
\end{abstract}
\maketitle
In their recent survey paper on enhancements for triangulated categories, Cano\-naco and Stellari pose the following
question:
\begin{question1} \cite[Question 3.13]{CanonacoStellari} Are there examples of triangulated categories, linear
  over a field $k$, with non-unique $k$-linear enhancement?
\end{question1}
Below we give such an example. In fact, we simply observe that the
topological graded field example (see \cite[\S2.1]{Schwede}) can be
made to work in the algebraic case\footnote{Another, more complicated,
  but also potentially more interesting, example is presented in~\cite{Kajiura}. The techniques in loc.\ cit.\ are based on an
  ingeneous direct manipulation of solutions to the Maurer-Cartan
  equation.  However the authors feel that as it stands, the arguments
  are not fully complete. In particular while it is possible to
  ``remove strict units'' \cite[\S4.4]{Kajiura}) from solutions to the
  Maurer-Cartan solution, the required $A_\infty$-isomorphism will in
  general be more complicated than the proof suggests.  This makes the
  verification of ``Condition (1)'' in loc.\ cit.\ more delicate (if at all possible)
  and therefore the same is true for the claim that the constructed
  functor is exact. We are currently discussing these points with the
  author.
}.
\begin{example1} 
Let $K=k(x_1,\cdots,x_{n+1})$  and $F=K[t,t^{-1}]$, where $n>0$ is even,
$t$ has cohomological degree $n$ and $K$ is concentrated in degree zero.  Since all homogeneous elements in $F$ are invertible, we call $F$ a graded\footnote{Throughout all graded notions are interpreted in the ``super'' sense. This affects for example the signs in the Hochschild complex.} field.
Let\footnote{\label{fn} We write $\HH^l(A)$ for the $l$'th graded Hochschild cohomology group of a graded $k$-algebra $A$. In particular $\HH^l(A)$ is a $\ZZ$-graded
  vector space. If we consider
  $A$ as a dg-algebra $A_{\mathrm dg}$ with trivial differential then $\HH^l(A_{\mathrm dg})=\prod_{i+j=l} \HH^i(A)^j$.}  $0\neq \eta\in \HH^{n+1}_k(K)$.
Then, by Lemma 1 below, $\tilde{\eta}=\eta\otimes d/dt$ is a non-zero element of
$\HH^{n+2}_k(F)^{-n}$. Since, also by Lemma 1, $\HH^{s}_k(F)=0$ for $s>n+2$, we may
construct a minimal $A_\infty$ structure
$(0,m_2,0,\ldots,0,m_{n+2},m_{n+3},\ldots)$ on $F$ such that the class
of $m_{n+2}$ is $\tilde{\eta}$ (see \cite[Lem.\ B.4.1]{Lefevre}). Let $F_\eta$
be the resulting $A_\infty$-algebra and let $f_1\in \Aut_k(F)$. One checks using \cite[Lem.\ B.4.2]{Lefevre}) or directly that $\tilde{\eta}\circ f_1$ is the first obstruction
against extending $f_1$ to an $A_\infty$-isomorphism $f:F\r F_\eta$. Since $\tilde{\eta}$ is non-trivial,
the same is true for $\tilde{\eta}\circ f_1$ and so $F$ and $F_\eta$ are not $A_\infty$-isomorphic.

As in \cite{Lurie1}, we see that the  triangulated category $\Perf(F_\eta)$ of right
perfect $F_\eta$-modules is equivalent, as a graded category, to the category of graded $F$-vector spaces
of finite rank. Since the latter category is semi-simple, it has only one triangulated
structure compatible with the graded structure.
 Hence
$\Perf(F)$ and $\Perf(F_\eta)$ are equivalent as trianguled
categories.

On the other hand, $\Perf(F)$ and $\Perf(F_\eta)$ have canonical
$A_\infty$-enhancement given by the $A_\infty$-categories of twisted complexes $\Tw(F)$ and $\Tw(F_\eta)$ (see \cite[Ch.\ 7]{Lefevre}).  We claim that $\Tw(F)$
and $\Tw(F_\eta)$ are not $A_\infty$-equivalent. Indeed, any  $A_\infty$-equivalence between them
would have to send the indecomposable (right) $F$-module $F_F$ to an object in $\Tw(F_\eta)$ which is $A_\infty$-isomorphic to $\Sigma^u(F_\eta)_{F_\eta}$ for some $u$. Hence
$\Tw(F)(F_F,F_F)\cong F$ and
$\Tw(F_\eta)(\Sigma^u(F_\eta)_{F_\eta},\Sigma^u(F_\eta)_{F_\eta})\cong F_\eta$ would have to be
$A_\infty$-isomorphic $A_\infty$-algebras (since they are both minimal). This is not the
case as we have established above.
\end{example1}
\begin{remark1} There is nothing special about the particular pair $(K,\eta)$ we have used. The chosen $(K,\eta)$ simply allows
for the most trivial argument for the existence of an $A_\infty$-structure on $F$ with the given $m_{n+2}$.
\end{remark1}
We have used the following basic lemma:
\begin{lemma1}
\label{lem:hochschild}
Let $F=K[t,t^{-1}]$ be as above.
The Hochschild cohomology of $F$ is given by
$
  \HH^\ast_k(F)\cong\HH^\ast_k(K)\otimes_k \HH_k^\ast(k[t,t^{-1}])
$.
Moreover $\HH^i(k[t,t^{-1}])=0$ for $i>1$ and the derivation $d/d t$ represents
a non-trivial element of $\HH^1(k[t,t^{-1}])$.
\end{lemma1}
\begin{proof} By the next lemma we only have to understand $\HH^\ast(k[t,t^{-1}])$. 
Since $t$ has even degree, $F$ is graded commutative and so the claim follows from the graded version of the HKR theorem. 
\end{proof}
\begin{lemma2}\hspace*{-0.25cm}\footnote{This lemma is stated in \cite[Thm 4.7]{BerghOppermann} without any hypotheses on~$A,B$. However the proof in loc.\ cit.
(essentially the proof we have given)
requires some kind of finiteness hypothesis which seems to have been inadvertibly ommitted. Indeed the result is
false if $A$, $B$ are fields of infinite transcendence degree over $k$.

Note that this lemma does not automatically imply that $\HH^\ast((A\otimes_k B)_{\mathrm{dg}})=\HH^\ast(A_{\mathrm{dg}})\otimes_k \HH^\ast(B_{\mathrm{dg}})$
because of
the fact that tensor products need not commute with products  (cfr Footnote \ref{fn}). We thank Zhengfang Wang for this observation.} Let $A,B$ be graded $k$-algebras such that the graded tensor product 
$B^e:=B\otimes_k B^\circ$ is noetherian. Then $\HH^\ast(A\otimes_k B)\cong \HH^\ast(A)\otimes_k \HH^\ast(B)$.
\end{lemma2}
\begin{proof} 
 Let
$Q^\bullet$ be a resolution of $B$ by finitely generated graded projective $B$-bimodules, and let $P^\bullet$ be
an arbitrary resolution of $A$ by graded projective $A$-bimodules. Then $P^\bullet\otimes_k Q^\bullet$
is a graded projective resolution of $A\otimes_k B$ and we have
\begin{align*}
\HH^\ast(A\otimes_k B)&=H^\ast(\Hom_{A^e\otimes_k B^e}(P^\bullet\otimes_k Q^\bullet,A\otimes_k B))\\
&\cong H^\ast(\Hom_{A^e}(P^\bullet,A)\otimes_k \Hom_{B^e}(Q^\bullet,B))\\
&\cong \HH^\ast(A)\otimes_k \HH^\ast(B).\qedhere
\end{align*}
\end{proof}
\def\cprime{$'$} \def\cprime{$'$} \def\cprime{$'$}
\providecommand{\bysame}{\leavevmode\hbox to3em{\hrulefill}\thinspace}
\providecommand{\MR}{\relax\ifhmode\unskip\space\fi MR }
\providecommand{\MRhref}[2]{%
  \href{http://www.ams.org/mathscinet-getitem?mr=#1}{#2}
}
\providecommand{\href}[2]{#2}

\end{document}